\documentclass[letterpaper, 10 pt, conference]{ieeeconf}  

\IEEEoverridecommandlockouts                 

\usepackage{amsmath}
\usepackage{amssymb}
\usepackage{graphicx}
\usepackage[colorlinks=true, allcolors=blue]{hyperref}
\usepackage{makecell}
\usepackage{multirow}

\usepackage{color}

\usepackage{graphicx}
\usepackage{subcaption}
\usepackage{booktabs}

\title{\LARGE \bf
On Tikhonov Regularization for Direct and Indirect \\ Data-Driven LQR Control
}

\author{Shuyuan Zhang$^{1}$, Zheming Wang$^{2}$, and Rapha\"{e}l M. Jungers$^{1}$
\thanks{RJ is a FNRS honorary Research Associate. This project has received funding from the \textit{European Research Council (ERC) under the European Union's Horizon 2020 research and innovation programme} under grant agreement No 864017 - L2C, from the Horizon Europe programme under grant agreement No101177842 - Unimaas, and from the ARC (French Community of Belgium) - project name: SIDDARTA. 
The work of Zheming Wang was supported in part by the National Natural Science Foundation of China under Grant 62303416 and W2421028.}
\thanks{$^{1}$S. Zhang and R. M. Jungers are with the ICTEAM Institute, UCLouvain, 4 Avenue Georges Lema\^itre, 1348 Louvain-la-Neuve, Belgium. {\tt\small shuyuan.zhang,raphael.jungers@uclouvain.be}} 
\thanks{$^{2}$Z. Wang is with the Department of Automation, Zhejiang University of Technology, Hangzhou 310023, China
{\tt\small wangzheming@zjut.edu.cn}}
}

\begin{document}

\maketitle
\thispagestyle{empty}
\pagestyle{empty}

\begin{abstract}
In recent years, the so-called `direct data-driven control' has been a topic of intense research, and it is expected that it will become prominent in future complex dynamical systems control. 
Within this framework, regularization not only implicitly enforces system identification, but also plays a crucial role in ensuring reliable closed-loop behavior.
To further enhance the performance of data-driven controllers, we propose a new regularization method for direct data-driven LQR control of unknown LTI systems, based on a regularized covariance parameterization.
Unlike existing data-driven techniques, the proposed method remains effective in handling ill-conditioned cases, such as when the data matrix has a large condition number.
Then, we demonstrate that our method is equivalent to the indirect certainty-equivalence LQR combined with Tikhonov regularization.
Furthermore, we extend our method to the design of controllers for unknown nonlinear systems using Koopman linear embedding. 
Finally, the simulation results validate the effectiveness and advantages of the proposed regularization method.


\end{abstract}

\section{INTRODUCTION}

Over the past two decades, data-driven methods have gained increasing popularity, particularly in cases where no available model or expert knowledge exists \cite{hou2013model,martin2023guarantees}. By learning from data, these methods open new avenues for developing intelligent solutions in modern control systems. Broadly, data-driven methods can be classified into two main types: indirect and direct.

\textit{Indirect methods} rely on system identification \cite{ljung1999system}. 
A model is first identified from data, and then analyzed in a second phase. In this phase, a controller is designed using traditional techniques, such as Linear Quadratic Regulator (LQR), treating the identified system from the first phase as the ground-truth \cite{mania2019certainty}. 
This approach relies on the so called `certainty-equivalence principle' \cite{astrom1994adaptive}, whereby the estimates are treated as if they were the true parameters, without accounting for their uncertainty. 
In the first phase, regularized identification methods \cite{pillonetto2022regularized} can be applied to regression problems to handle uncertainties and facilitate numerical computation.

\textit{Direct methods} bypass the system identification step \cite{maupong2017data}, relying solely on measured input-output data from the controlled plant for controller design. In \cite{willems2005note}, Willems \textit{et al.} presented a foundational result, known as Willems' fundamental lemma, which implies that the behavior of an LTI system can be characterized by the column span of a data Hankel matrix.
Inspired by this lemma, De Persis and Tesi \cite{de2020formulas} proposed a data-dependent representation of the closed-loop behavior of an LTI system for LQR design by solving a positive semi-definite program.
Furthermore, they developed a regularizer that promotes robust stability when the data used for controller design are affected by noise \cite{de2021low}.
Subsequent works \cite{dorfler2023bridging,dorfler2022role,dorfler2023certainty,zhao2025regularization} investigated different regularizers aimed at enhancing system performance.
In \cite{dorfler2023certainty}, D{\"o}rfler \textit{et al.} presented a regularized direct data-driven LQR formulation, which has been proven to be equivalent to the indirect certainty-equivalence LQR formulation with ordinary least-squares identification (see Corollary III.2 and Theorem III.3 in \cite{dorfler2023certainty}). This method allows a trade-off between performance and robustness by blending different regularizers that promote either certainty equivalence or robust stability.
Recently, Zhao \textit{et al.} \cite{zhao2025regularization} proposed a regularization method for direct data-driven LQR design using a covariance parameterization. A notable advantage of this parameterization is that the dimension of the resulting formulation is independent of the data length. This approach has been applied to autonomous bicycle control \cite{persson2025adaptive} and power converters \cite{zhao2024direct}.

While the parameterization and regularization approaches shed light on direct data-driven LQR design, they rely on the assumption that the available data is sufficiently rich. 
This assumption guarantees that the data matrix $D_{0}D_{0}^{\top}$ is positive definite, as shown in \eqref{V}. 
However, its condition number may still be large, and inverting it may lead to large numerical errors.
To address this issue, inspired by Tikhonov regularization \cite{hoerl1970ridge}, this paper proposes a regularized covariance parameterization for direct data-driven LQR design of unknown LTI systems corrupted by process noise.
Unlike the approaches in the above literature, the proposed parameterization remains valid even when $D_{0}D_{0}^{\top}$
has a large condition number or when $D_{0}$ is rank-deficient. 
Moreover, we show that the proposed method is equivalent to the indirect certainty-equivalence LQR approach combined with Tikhonov regularization (i.e., regularized least-squares identification). 
Notably, the dimension of the resulting LQR formulation remains independent of the data length.
Compared with the regularized direct data-driven control method in \cite{zhao2025regularization}, our approach achieves superior control performance through tuning of the regularization coefficients (see Examples 1 and 2 for details).
Finally, we extend the proposed method to controller design for unknown nonlinear systems using Koopman linear embedding. 


The rest of this paper is organized as follows. Section~\ref{2} reviews the model-based and data-driven LQR formulations. Section~\ref{3} presents a regularization method for direct data-driven LQR design of LTI systems, based on the proposed regularized covariance parameterization. In Section~\ref{NON}, the proposed method is extended to handle nonlinear systems. 
Section \ref{4} validates the theoretical results through simulations. Conclusions are drawn in Section \ref{5}.

\section{PRELIMINARIES}
\label{2}

In this section, we review the model-based LQR, the indirect certainty-equivalence LQR \cite{dorfler2022role,dorfler2023certainty}, and the direct data-driven LQR design via covariance parameterization along with its regularized variant \cite{zhao2025regularization,zhao2025data}.

\subsection{Model-Based LQR}
Consider a discrete-time LTI system:
\begin{equation} \label{LTI}
    x(k+1) = Ax(k) + Bu(k) + w(k)
\end{equation}
with $k \in \mathbb{N}$, where $x(k) \in \mathbb{R}^n$, $u(k) \in \mathbb{R}^m$ and $w(k) \in \mathbb{R}^n$ are the state, input and noise at time $k$, respectively. The pair $(A,B)$ is unknown but assumed to be stabilizable. 

{\bf Assumption 1.} The initial state and the process noise are independent random variables with distributions satisfying $x(0) \overset{\text{i.i.d.}}{\sim} \mathcal{N}(0, \sigma_{x}^{2} I_{n})$, and $w(k) \overset{\text{i.i.d.}}{\sim} \mathcal{N}(0, \sigma_{w}^{2} I_{n})$ for $k \in \mathbb{N}$, where $\sigma_{x}>0$ and $\sigma_{w}>0$.

The infinite-horizon LQR problem aims to find a state feedback gain $K$ that minimizes the time-average cost, i.e.,
\begin{align} 
& \min_K J(K) := \lim_{T \to \infty} \mathbb{E} \left[\frac{1}{T} \sum_{k=0}^{T-1} C(k) \right] \nonumber\\
& \text{s.t.} \ \eqref{LTI} \ \text{and} 
\ u(k) = K x(k),     
\end{align}
where $C(k) = x(k)^\top Q x(k) + u(k)^\top R u(k)$ with $Q > 0$ and $R > 0$, and $\mathbb{E}$ denotes the expectation over the randomness from the initial state and the process noise.

When $A + BK$ is stable, it holds that \cite{anderson2007optimal}
\begin{equation} \label{Cost}
    J(K) = \text{Tr}(QP + K^{\top} R K P),
\end{equation}
where $P$ is the unique positive definite solution to the discrete-time Lyapunov equation:
\begin{equation} \label{Riccati}
    P = I_{n} + (A + BK) P (A + BK)^{\top}.
\end{equation}

The optimal LQR gain $K^{*} := \arg \min_{K} J(K)$ of \eqref{Cost} and \eqref{Riccati} can be found by the celebrated Riccati equation solution \cite{anderson2007optimal}. Alternatively, as shown in Section 2 of \cite{4792678}, $K^{*}$ can be determined by solving the following optimization problem:
\begin{align} \label{gruond_K}
& \min_{P\geq I_{n}, K} \text{Tr}(QP + K^{\top} R K P) \nonumber\\
& \text{s.t.} \ (A + BK) P (A + BK)^{\top} - P +  I_{n} \leq 0.  
\end{align}
We refer to \eqref{gruond_K} as the model-based LQR problem. 

However, when $(A, B)$ is unknown, a direct solution to \eqref{gruond_K} is not feasible. In the following,  we recapitulate several data-driven control methods for computing $K$ without explicit knowledge of $(A, B)$.

\subsection{Indirect Certainty-Equivalence LQR}
\label{sec2-indirect}

A classical approach to data-driven LQR design is based on the certainty-equivalence principle \cite{mania2019certainty}: it first estimates a model from data, and then solves the LQR problem regarding the identified model as the ground-truth.

We assume that the data are collected offline by exciting system \eqref{LTI} with a random input $u(k)$ for $k \in \mathbb{N}$, 
where the input is independent of both the process noise and the initial condition specified in Assumption 1.
Without loss of generality, we consider a $T$-step time series of states, inputs, noises, and successor states:
\begin{align*}
X_{0} &:= \begin{bmatrix} 
x(0) & x(1) & \cdots & x(T-1) 
\end{bmatrix} \in \mathbb{R}^{n \times T},  \\
U_{0} &:= \begin{bmatrix} 
u(0) & u(1) & \cdots & u(T-1) 
\end{bmatrix} \in \mathbb{R}^{m \times T},  \\
W_{0} &:= \begin{bmatrix} 
w(0) & w(1) & \cdots & w(T-1) 
\end{bmatrix} \in \mathbb{R}^{n \times T},  \\
X_{1} &:= \begin{bmatrix} 
x(1) & x(2) & \cdots & x(T) 
\end{bmatrix} \in \mathbb{R}^{n \times T}.
\end{align*}
Thus, we have a data-based representation of system \eqref{LTI} as follows: 
\begin{equation} \label{data LTI}
    X_{1} = \begin{bmatrix} B & A \end{bmatrix} 
    \begin{bmatrix} U_{0} \\ X_{0} \end{bmatrix} + W_{0}.
\end{equation}
Let us define
\[
    D_{0} := \begin{bmatrix} U_{0} \\ X_{0} \end{bmatrix}.
\]

{\bf Assumption 2.} The collected data is sufficiently rich, that is,
\begin{equation} \label{rank}
    \text{rank}(D_{0}) = n + m.
\end{equation}
This assumption is a necessary condition for the existing data-driven LQR design approach \cite{van2020data}.

Based on \eqref{data LTI} and \eqref{rank}, the least-squares estimator $(\hat{A}, \hat{B})$ is given by
\begin{equation} \label{LS}
    [\hat{B}, \hat{A}] = \arg \min_{B,A} \| X_{1}-[B, A] D_{0} \|_{F} = X_{1} D_{0}^{\dagger},
\end{equation}
where $D_{0}^{\dagger}=D_{0}^{\top}(D_{0}D_{0}^{\top})^{-1}$. 

As presented in \cite{dorfler2023certainty}, the system $(A,B)$ is replaced with its estimate $(\hat{A}, \hat{B})$ in \eqref{gruond_K}, and the LQR problem can be reformulated as
\begin{align} \label{indirect}
& \min_{P\geq I_{n}, K} \text{Tr}(QP + K^{\top} R K P) \nonumber\\
& \text{s.t.} \ (\hat{A} + \hat{B}K) P (\hat{A} + \hat{B}K)^{\top} - P +  I_{n} \leq 0.
\end{align}
We refer to \eqref{indirect} as the indirect certainty-equivalence LQR problem.

\subsection{Direct Data-Driven LQR with Covariance Parameterization}
\label{2-3}

The goal of direct data-driven LQR design is to compute the feedback gain $K$ without performing explicit system identification \eqref{LS}.
This subsection reviews recent advances in direct data-driven control methods based on parameterization and regularization \cite{zhao2025regularization,zhao2025data}.

Define $\Phi := D_{0} D_{0}^{\top} / T$. 
Under Assumption 2, $\Phi$ is positive definite, and for every $K$ there exists a unique $V \in \mathbb{R}^{(n+m)\times n}$ such that 
\begin{align} \label{V}
\begin{bmatrix}
K \\
I_n
\end{bmatrix}
= \Phi V.  
\end{align}
We refer to \eqref{V} as the covariance parameterization \cite[(15)]{zhao2025data}. Compared to the data-based parameterization
\begin{align} \label{G}
\begin{bmatrix}
K \\
I_n
\end{bmatrix}
= D_{0} G, \quad G \in \mathbb{R}^{T \times n},
\end{align}
proposed in \cite[(9)]{dorfler2023certainty}, the dimension of \( V \) is independent of the data length.

Define the sample covariance matrices $\bar{X}_{0} = X_{0} D_{0}^{\top}/T$, $\bar{U}_{0} = U_{0} D_{0}^{\top}/T$, $\bar{W}_{0} = W_{0} D_{0}^{\top}/T$, and $\bar{X}_{1} = X_{1} D_{0}^{\top}/T$.  
We can formulate the direct data-driven LQR problem as
\begin{align} \label{direct}
& \min_{P\geq I_{n}, V} \text{Tr}(QP + (\bar{U}_{0} V)^{\top} R (\bar{U}_{0} V) P) \nonumber\\
& \text{s.t.} \  (\bar{X}_{1} V) P (\bar{X}_{1} V)^{\top} - P +  I_{n} \leq 0, \nonumber\\
& \quad \ \ \bar{X}_{0} V = I_{n}
\end{align}
with the gain matrix $K = \bar{U}_{0} V$. 
We refer to \eqref{direct} as the direct data-driven LQR problem with covariance parameterization (see (18) in \cite{zhao2025data}). 
As shown in \cite[Lemma 1]{zhao2025data}, \eqref{direct} is equivalent to the indirect certainty-equivalence LQR \eqref{indirect} in the sense that their solutions coincide.

Furthermore, Zhao \textit{et al.} \cite[(22)]{zhao2025regularization} presented a regularization formulation for the direct data-driven LQR problem \eqref{direct} to enhance robust closed-loop stability as follows:
\begin{align} \label{regularized direct}
& \min_{P\geq I_{n}, V} \text{Tr}(QP + (\bar{U}_{0} V)^{\top} R (\bar{U}_{0} V) P) + \lambda \Omega(V)\nonumber\\
& \text{s.t.} \  (\bar{X}_{1} V) P (\bar{X}_{1} V)^{\top} - P +  I_{n} \leq 0, \nonumber\\
& \quad \ \ \bar{X}_{0} V = I_{n}
\end{align}
with the regularization coefficient $\lambda \geq 0$, where $\Omega(V) := \text{Tr}(V P V^{\top} \Phi)$. 
We refer to \eqref{regularized direct} as the regularized direct data-driven LQR problem. 
As illustrated in \cite[Theorem 1]{zhao2025regularization}, the role of the regularizer $\Omega(V)$ is equivalent to mixing the certainty-equivalence promoting regularizer $\|\Pi G\|$ (with $\Pi := I_{T} - D_{0}^{\dagger}D_{0}$) from \cite{dorfler2023certainty} and the robustness-promoting regularizer $\text{Tr}(G P G^{\top})$ from \cite{de2021low}, where $G$ satisfies \eqref{G}. 

Various regularizers have been employed in direct data-driven LQR design to promote desired performance \cite{de2021low,dorfler2023certainty,zhao2025regularization}. 
In the next section, we introduce a new regularizer for direct data-driven LQR design, building upon the recent results in \cite{zhao2025regularization}.

\section{MAIN RESULTS}
\label{3}
In this section, we propose a new regularization method for direct data-driven LQR design and demonstrate its equivalence to the indirect certainty-equivalence LQR combined with Tikhonov regularization.

\subsection{Regularized Covariance Parameterization}
\label{3-1}

In this subsection, we present a direct data-driven LQR design method using the regularized covariance parameterization to compute the controller. 

The covariance parameterization \eqref{V} relies on Assumption 2, which prevents $D_0 D_0^{\top}$ from having exactly zero eigenvalues. However, the condition number of $D_0 D_0^\top$, denoted by $\text{cond}(D_0 D_0^\top)$, can still be large, and inverting $D_0 D_0^\top$ may lead to large numerical errors.

Inspired by the Tikhonov regularization technique \cite{hoerl1970ridge} used in system identification, we add a \emph{ridge} term to the diagonal, i.e., $D_{0}D_{0}^{\top}+ \gamma I_{n+m}$ with $\gamma > 0$, which shifts all eigenvalues of $D_0 D_0^{\top}$ by $\gamma$ and thereby mitigates ill-conditioning.
To this end, we define $\Psi := (D_{0}D_{0}^{\top}+\gamma I_{n+m})/T$. 
Since $\Psi$ is positive definite even when $\text{rank}(D_{0}) < n + m$ (i.e., Assumption 2 is not satisfied), it follows that for every $K$ there is a unique $\Xi \in \mathbb{R}^{(n+m)\times n}$ such that
\begin{align} \label{RCP}
\begin{bmatrix}
K \\
I_n
\end{bmatrix}
= \Psi \Xi.  
\end{align}
We refer to \eqref{RCP} as the regularized covariance parameterization.

Define $\Psi:= \begin{bmatrix}
\Psi_{1} \\
\Psi_{2}
\end{bmatrix}$, where $\Psi_{1} \in \mathbb{R}^{m \times (n+m)}$ and $\Psi_{2} \in \mathbb{R}^{n \times (n+m)}$. 
Similar to \eqref{direct}, we formulate the following direct data-driven LQR problem:
\begin{align} \label{ridge direct}
& \min_{P\geq I_{n}, \Xi} \text{Tr}(QP + (\Psi_{1} \Xi)^{\top} R (\Psi_{1} \Xi) P) \nonumber\\
& \text{s.t.} \  (\bar{X}_{1} \Xi) P (\bar{X}_{1} \Xi)^{\top} - P +  I_{n} \leq 0, \nonumber\\
& \quad \ \ \Psi_{2} \Xi = I_{n}
\end{align}
with the gain matrix $K = \Psi_{1} \Xi$. We refer to \eqref{ridge direct} as the direct data-driven LQR problem with regularized covariance parameterization.

{\bf Remark 1.} The dimensions of all decision variables in \eqref{ridge direct} depend only on the system parameters $n$ and $m$, and not on the data length $T$. This method is well-suited for larger datasets. 
By appropriately tuning \( \gamma \), our method can reliably produce a stabilizing controller (see Example~1 in Section~\ref{4}). In particular, when \( \gamma = 0 \), we have \( \Psi_{1} = \bar{U}_{0} \) and \( \Psi_{2} = \bar{X}_{0} \), and \eqref{ridge direct} reduces to \eqref{direct}.

{\bf Remark 2.} Although both regularized direct data-driven formulations \eqref{regularized direct} and \eqref{ridge direct} simultaneously address identification and control objectives, our method remains effective in handling ill-conditioned cases, e.g., when $\text{cond}(D_0 D_0^\top)$ is large.
In addition, we find that the control performance can be enhanced by mixing the two regularizers in \eqref{regularized direct} and \eqref{ridge direct}, which is demonstrated in Example 1 of Section \ref{4}.

\subsection{Equivalence to Indirect Certainty-Equivalence LQR}
\label{3-2}

This subsection establishes the equivalence between the proposed \eqref{ridge direct} and the indirect certainty-equivalence LQR combined with Tikhonov regularization.

In system identification, a classical approach is Tikhonov regularization method \cite{hastie2020ridge}, which penalizes large values of $[B,A]$ in the following optimization problem:
\begin{align} \label{ReLS}
[\hat{B}, \hat{A}] 
&= \arg \min_{B,A} \| X_{1}-[B, A] D_{0} \|_{2}^{2} + \gamma ||[B, A]||_{2}^{2} \nonumber\\
&= X_{1} D_{0}^{\top}(D_{0}D_{0}^{\top}+\gamma I_{n+m})^{-1}. 
\end{align} 
This is the ridge regression solution proposed by Hoerl and Kennard \cite{hoerl1970ridge} in the 1970s. 
The regularized least-squares estimator \eqref{ReLS} allows to handle rank-deficient data matrices, i.e., when $\text{rank}(D_{0}) < n + m$. The positive definiteness of $D_{0}D_{0}^{\top}+ \gamma I_{n+m}$ can be ensured by tuning $\gamma$. Consequently, the estimator \eqref{ReLS} remains applicable even when Assumption 2 is not satisfied.

Similar to \eqref{indirect}, one can replace $(A,B)$ with its estimate $(\hat{A}, \hat{B})$ in \eqref{ReLS}, which yields the following LQR problem:
\begin{align} \label{Ridge indirect}
& \min_{P\geq I_{n}, K} \text{Tr}(QP + K^{\top} R K P) \nonumber\\
& \text{s.t.} \ (\hat{A} + \hat{B}K) P (\hat{A} + \hat{B}K)^{\top} - P +  I_{n} \leq 0.
\end{align}
We refer to \eqref{Ridge indirect} as the regularized indirect certainty-equivalence LQR problem.

Next, we establish the equivalence between the formulations \eqref{ridge direct} and \eqref{Ridge indirect}.

{\bf Theorem 1.} Consider the regularized direct and indirect data-driven LQR formulations \eqref{ridge direct} and \eqref{Ridge indirect}, respectively. These two formulations are equivalent in the sense that their solutions coincide.

{\bf Proof.} Due to the positive definiteness of $\Psi$, it follows from \eqref{RCP} that $\Xi$ is uniquely determined as $\Xi = \Psi^{-1} \begin{bmatrix}
K \\
I_n
\end{bmatrix}$. Then, the LQR problem \eqref{ridge direct} can be reformulated as
\begin{align*} 
& \min_{P\geq I_{n}, K} \text{Tr}(QP + K^{\top} R K P) \nonumber\\
& \text{s.t.} \ X_{1} D_{0}^{\top} (D_{0}D_{0}^{\top}+\gamma I_{n+m})^{-1} \begin{bmatrix}
K \\
I_n
\end{bmatrix} P \begin{bmatrix}
K \\
I_n
\end{bmatrix}^{\top} \nonumber\\
& \quad \quad \ \ (X_{1} D_{0}^{\top} (D_{0}D_{0}^{\top}+\gamma I_{n+m})^{-1})^{\top} - P +  I_{n} \leq 0. 
\end{align*}
By the definition $[\hat{B}, \hat{A}] 
= X_{1} D_{0}^{\top}(D_{0}D_{0}^{\top}+\gamma I_{n+m})^{-1}$ in \eqref{ReLS}, this problem is
exactly the regularized indirect certainty-equivalence LQR problem \eqref{Ridge indirect}.
\hfill{$\blacksquare$}

Furthermore, we show that problem \eqref{ridge direct} can be cast as a convex program.

{\bf Proposition 1.}
Suppose that \eqref{ridge direct} is feasible. Then, the optimal gain for \eqref{ridge direct} can be computed as $ K = \Psi_{1} Y P^{-1}$, where \( Y \in \mathbb{R}^{(n+m) \times n} \) and \( P \in \mathbb{R}^{n \times n} \) are solutions to
\begin{align} \label{TH2}
&\underset{P,Y,L}{\min} 
\operatorname{Tr}(Q P) + \operatorname{Tr}(RL) \nonumber\\
&\text{s.t.} 
\left\{ 
\begin{aligned}
&\Psi_{2} Y = P, \\
&\begin{bmatrix} 
P - I_n & \bar{X}_{1} Y \\ 
Y^{\top} \bar{X}_{1}^{\top} & P 
\end{bmatrix} \geq 0, \\
&\begin{bmatrix} 
L &  \Psi_{1} Y \\ 
Y^{\top} \Psi_{1}^{\top} & P 
\end{bmatrix} \geq 0, 
\end{aligned}
\right.
\end{align}
where \( L \in \mathbb{R}^{m \times m} \).

{\bf Proof.} We first show that the solution to \eqref{ridge direct} can be computed as $ K = \Psi_{1} Y P^{-1}$, where $Y$ and $P$ are solutions to
\begin{align} \label{TH2-1}
&\underset{P,Y,L}{\min} 
\operatorname{Tr}(Q P) + \operatorname{Tr}(RL) \nonumber\\
&\text{s.t.} 
\left\{ 
\begin{aligned}
&\Psi_{2} Y = P, \\
&P \geq I_{n}, \\
&\bar{X}_{1} Y P^{-1} Y^{\top} \bar{X}_{1}^{\top} - P + I_{n} \leq 0, \\
&L - \Psi_{1} Y P^{-1} Y^{\top} \Psi_{1}^{\top} \geq 0.
\end{aligned}
\right.
\end{align}
By the variable substitution $\Xi = Y P^{-1}$, the first three constraints of \eqref{TH2-1} are equivalent to those of \eqref{ridge direct}. Hence, the feasibility of \eqref{TH2-1} implies that of \eqref{ridge direct}. 
The converse also holds by constructing an $L$ that satisfies the fourth constraint of \eqref{TH2-1}.
Let \( (P,Y,L) \) be an optimal solution to \eqref{TH2-1} and \( (P, \Xi) \) an optimal solution to \eqref{ridge direct}. Since \( R \) is positive definite, the optimal solution satisfies \( L = \Psi_{1} Y P^{-1} Y^{\top} \Psi_{1}^{\top} \). Therefore, the optima of \eqref{TH2-1} and \eqref{ridge direct} coincide.

Based on the Schur complement \cite{gallier2010notes}, \eqref{TH2-1} can be reformulated as \eqref{TH2}. Substituting \( \Xi = Y P^{-1} \) and \( K = \Psi_{1} \Xi \) then completes the proof.
\hfill{$\blacksquare$}

\section{EXTENSIONS to NONLINEAR SYSTEMS}
\label{NON}

In this section, we extend our method to the design of controllers for unknown nonlinear systems.

Consider a discrete-time nonlinear system
\begin{equation} \label{nonlinear}
    x(k+1) = f(x(k), u(k)), 
\end{equation}
with $k \in \mathbb{N}$, where $x(k) \in \mathbb{R}^n$ and $u(k) \in \mathbb{R}^m$ are the state and input at time $k$, respectively, and $f(\cdot)$ is the unknown nonlinear dynamics. 

One key idea of the Koopman operator is to lift the state $x(k)$ of the nonlinear system \eqref{nonlinear} to a higher-dimensional space via a set of lifting functions (often referred to as \emph{observables}) \cite{korda2018linear},
where the evolution of these observables becomes linear.

{\bf Assumption 3.} The nonlinear system \eqref{nonlinear} is assumed to admit a Koopman linear embedding. 
That is, there exists a set of linearly independent lifting functions 
$\phi_1(\cdot), \ldots, \phi_{n_z}(\cdot): \mathbb{R}^n \to \mathbb{R}$ 
such that the lifted state
\begin{equation}
    \Theta(x(k)) := [\phi_{1}(x(k)), \ldots, \phi_{n_{z}}(x(k))] \in \mathbb{R}^{n_{z}}
    \label{lifted_state}
\end{equation}
evolves linearly along all trajectories of \eqref{nonlinear}.

Under Assumption 3, the lifted state $z(k) = \Theta(x(k)) \in \mathbb{R}^{n_{z}}$ with $n_{z} > n$ satisfies
\begin{equation} \label{koopman_linear_model}
    z(k+1) = A z(k) + B u(k),     
\end{equation}
with matrices $A $ and $B$ of appropriate dimensions.


For system \eqref{nonlinear}, we collect a time series of length $T$ of states, inputs, and successor states, denoted by $X_{0}$, $U_{0}$, and $X_{1}$, as in Section \ref{sec2-indirect}. 
With the lifting functions \eqref{lifted_state}, we compute the lifted states as
\begin{align*}
    Z_{0} &:= \begin{bmatrix} 
    \Theta(x(0)) & \cdots & \Theta(x(T-1))
    \end{bmatrix}  \in \mathbb{R}^{n_{z} \times T}, \\
    Z_{1} &:= \begin{bmatrix} 
    \Theta(x(1)) & \cdots & \Theta(x(T)) \end{bmatrix}  \in \mathbb{R}^{n_{z} \times T}.
\end{align*}
Let 
\[
    D_{0} := \begin{bmatrix} U_{0} \\ Z_{0} \end{bmatrix}, 
    \quad 
    \Psi := \frac{1}{T}\big(D_{0}D_{0}^{\top} + \gamma I_{n+m}\big).
\]
Define $\bar{Z}_{1} = Z_{1} D_{0}^{\top}/T$. 
By solving \eqref{ridge direct} with $\bar{Z}_{1}$ in place of $\bar{X}_{1}$, the feedback gain $K$ is obtained. This yields the controller $u(k) = K z(k)$ for the linear system \eqref{koopman_linear_model}, and consequently applies to the original nonlinear system \eqref{nonlinear}.




\section{SIMULATIONS}
\label{4}

This section presents numerical simulations to validate the effectiveness and demonstrate the advantages of the proposed theoretical results. 

{\bf Example 1.} Consider the following linear system, where the unknown matrices $A$ and $B$ are chosen as:
\begin{align}  \label{AB}
A =
\begin{bmatrix}
0.99 & 0.01 & 0.02 & 0 \\
0.02 & 0.98 & 0.01 & 0 \\
0.01 & 0.03 & 0.97 & 0.01 \\
0    & 0.01 & 0.02 & 0.95
\end{bmatrix}, \quad
B =
\begin{bmatrix}
1 \\[2pt]
0 \\[2pt]
0 \\[2pt]
0
\end{bmatrix}.
\end{align}

This system has $n = 4$, and $m = 1$, with an unstable open-loop behavior.
Let $Q = I_{n}$ and $R = 10^{-3} I_{m}$.
We perform 100 independent trials. In each trial, we run an experiment on the system with initial state $x(0) \sim \mathcal{N}(0, I_{n})$, input $u(k) \sim \mathcal{N}(0,1)$, and noise $w(k) \sim \mathcal{N}(0, \sigma_{w}^{2} I_{n})$ over a time horizon $T$. 
Let $K^{(i)}$ denote the controller obtained in the $i$-th trial by solving either \eqref{regularized direct} (referred to as Method I) or \eqref{ridge direct} (referred to as Method II), using MOSEK \cite{aps2019mosek}. 
We define the optimality gap as
\begin{equation}
E^{(i)} = \frac{J(K^{(i)}) - J(K^{*})}{J(K^{*})},
\label{eq:E}
\end{equation}
where $J(\cdot)$ is computed from \eqref{Cost}, and $K^{*}$ is obtained by solving \eqref{gruond_K}.
We denote by $S$ the percentage of trials yielding a stabilizing controller (with $S_{I}$ and $S_{II}$ corresponding to the two methods, respectively), and by $M$ the median value of $E^{(i)}$ over 100 independent trials (with $M_{I}$ and $M_{II}$).

\begin{figure}[!ht]
    \centering
    \begin{subfigure}{0.22\textwidth}
        \centering
        \includegraphics[width=\linewidth]{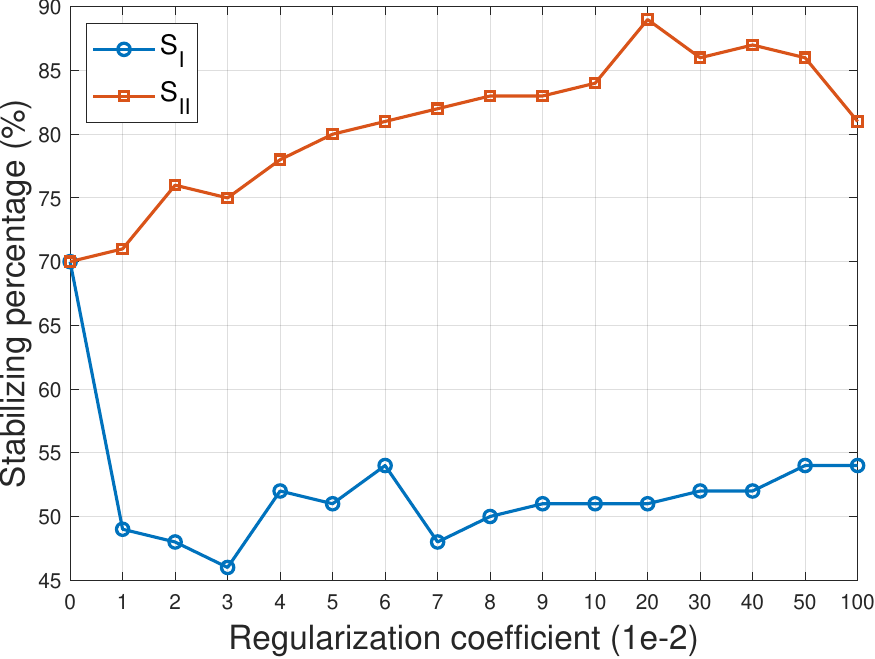}
        \caption{}
    \end{subfigure}
    \quad
    \begin{subfigure}{0.22\textwidth}
        \centering
        \includegraphics[width=\linewidth]{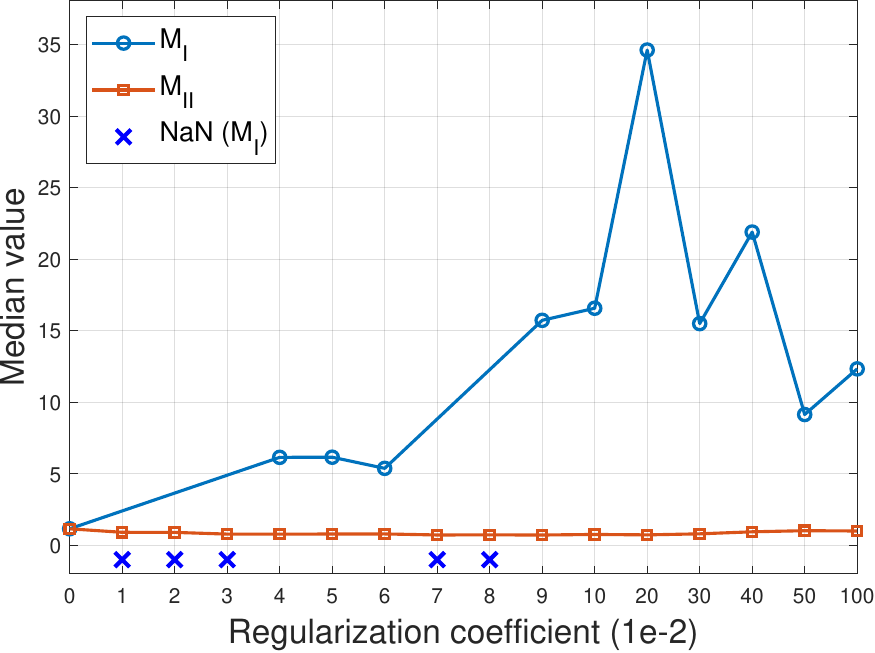}
        \caption{}
    \end{subfigure}

\caption{Comparison of Method~I and Method~II in terms of (a) stabilizing percentage $S$ and (b) median value $M$ as functions of regularization coefficients. Blue and orange lines correspond to $\lambda$ and $\gamma$, respectively.}

\label{EX1}
\end{figure}

Let $T = 10$ and $\sigma_{w} = 0.1$, and choose the regularization coefficients $\lambda$ and $\gamma$ from the set 
$[0,1,2,3,4,5,$ $6,7,8,9,10,20,30,40,50,100] \times 10^{-2}$. Fig.~\ref{EX1} shows the evolution of the stabilizing percentage $S$ and the median value $M$ with respect to the regularization coefficients, respectively.
Here, NaN indicates that no median value is reported because $S_{I} < 50$ for certain values of $\lambda$. 
We can observe that our Method II achieves better control performance than Method I, as indicated by a higher stabilizing percentage $S$ and a smaller median value $M$ through tuning of $\gamma$ (with $\max{(S_{II})} = 89$ and $\min{(M_{II})} = 0.7224$ across all $\gamma$).
The reason is that $\text{cond}(D_0 D_0^{\top}) \gg \text{cond}(D_0 D_0^{\top} + \gamma I_{n+m})$, with $\gamma$ improving the condition number of $D_0 D_0^{\top}$, whereas $\lambda$ does not. 
In addition, by combining the two regularizers, i.e., using 
$\text{Tr}(QP + (\Psi_{1} \Xi)^{\top} R (\Psi_{1} \Xi) P) + \lambda \Omega(\Xi)$ 
as the objective function in \eqref{ridge direct}, we find that the best performance is achieved $S = 93$ at $\lambda = 0.2$ and $\gamma = 0.3$, and $M = 0.6902$ at $\lambda = 0.03$ and $\gamma = 0.04$.
However, these regularizers cannot guarantee achieving the best performance in both $M$ and $S$ simultaneously.

\begin{figure}[!ht]
    \centering
    \begin{subfigure}{0.22\textwidth}
        \centering
        \includegraphics[width=\linewidth]{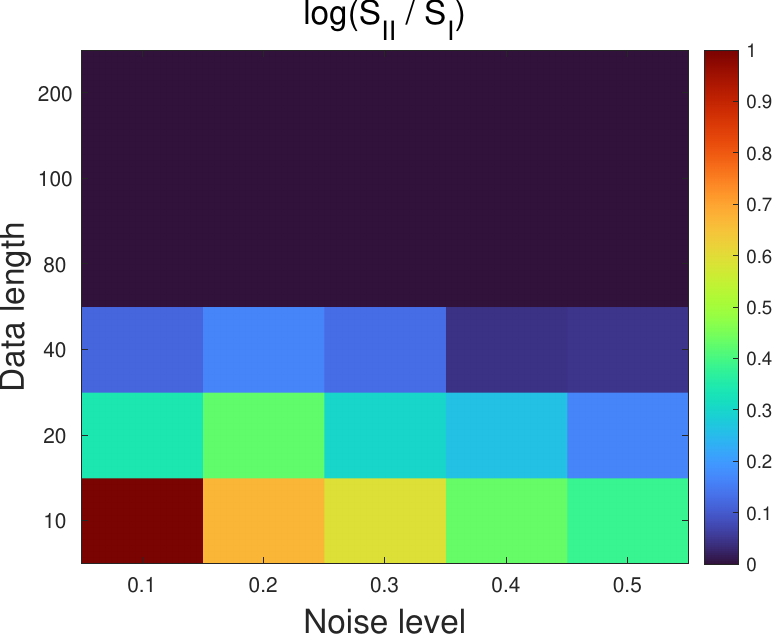}
        \caption{}
    \end{subfigure}
    \quad
    \begin{subfigure}{0.22\textwidth}
        \centering
        \includegraphics[width=\linewidth]{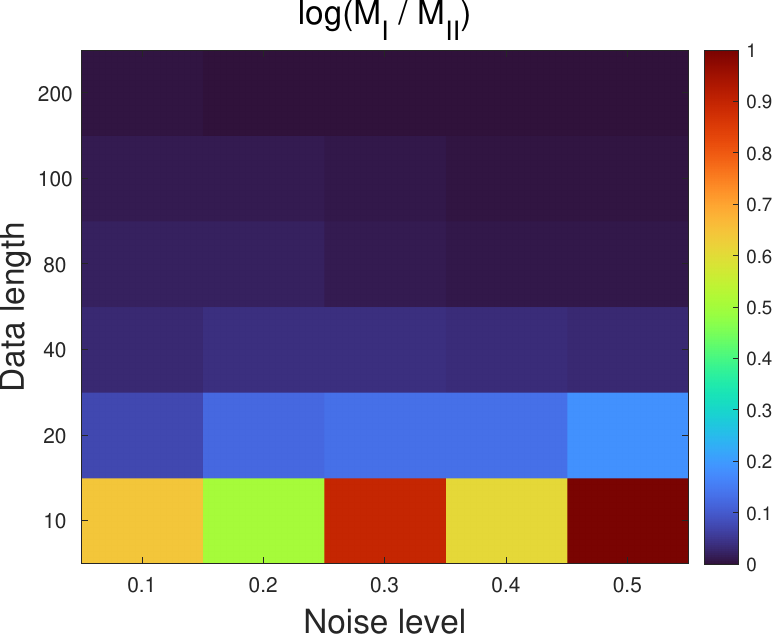}
        \caption{}
    \end{subfigure}

\caption{(a) $\log(S_{II}/S_{I})$ as a function of $T$ and $\sigma_{w}$; (b) $\log(M_{I}/M_{II})$ as a function of $T$ and $\sigma_{w}$.}
\label{EX1-heat}
\end{figure}

We now increase the data length and noise level, and choose $T$ and $\sigma_{w}$ from the following sets: $T = [10, 20, $ $40, 80, 100, 200]$ and $\sigma_{w} = [0.1, 0.2, 0.3, 0.4, 0.5]$, respectively.
For each pair $(T, \sigma_{w})$, the colors in Fig.~\ref{EX1-heat} represent the logarithmic values of $S_{II}/S_{I}$ and $M_{I}/M_{II}$. 
For each $\gamma$, let $S_{II}^{\gamma}$ and $M_{II}^{\gamma}$ denote the stabilizing percentage and the median value of $J(K^{(i)})$ over 100 independent trials, respectively. 
Then, $S_{II}$ and $M_{II}$ represent the maximum and minimum of all $S_{II}^{\gamma}$ and $M_{II}^{\gamma}$, respectively.
The quantities $S_{I}$ and $M_{I}$ are defined in the same way as above, but with respect to $\lambda$.

As shown in Fig.~\ref{EX1-heat}, we have $\log(S_{II}/S_{I}) \geq 0$ and $\log(M_{I}/M_{II}) > 0$ for all $(T, \sigma_{w})$, indicating that Method~II consistently outperforms Method~I (For visualization, all values of $\log(S_{II}/S_{I})$ and $\log(M_{I}/M_{II})$ are normalized to the interval $[0,1]$).
As the data length $T$ increases, the colors become darker (blue), indicating that the ratios $S_{II}/S_{I}$ and $M_{I}/M_{II}$ decrease as the condition number of $D_0 D_0^{\top}$ improves, and both methods tend to exhibit comparable performance.
When $T = 200$, both methods yield nearly identical results; that is, $S_{II}$ and $S_{I}$ are 100\%, and $M_{I}$ and $M_{II}$ are close to the ground-truth cost $J(K^{*})$.

Note that $T$ should not be chosen too large, since for unstable or stable systems the sampled trajectories will diverge to infinity or decay to zero, respectively, as 
$T$ increases significantly. In either case, the quality of the collected data deteriorates, and the resulting data matrix $D_0 D_0^{\top}$ may become ill-conditioned or numerically singular, making it infeasible to compute a stabilizing controller. 
Moreover, when $T$ is very small, Method I is not recommended, as its stabilizing percentage $S$ is too low, whereas Method II remains effective.
For instance, when $T = 5$ and $\sigma_{w} = 0.1$, $S_{I}$ is around $20\%$, while $S_{II}$ exceeds $70\%$.

\textbf{Example 2}. Consider 1000 random systems with $n = 3$ and $m=1$, where the entries of $A$ and $B$ are drawn from $\mathcal{N}(0,1)$. The data length is set to $T = 10$. 

For each system, the regularization coefficient is randomly selected from the interval $(0,1)$. We then perform 50 independent trials, where $S$ denotes the percentage of trials that yield a stabilizing controller. In each trial, we run an experiment on the system with initial state $x(0) \sim \mathcal{N}(0, I_{n})$, input $u(k) \sim \mathcal{N}(0, 1)$, and noise $w(k) \sim \mathcal{N}(0, \sigma_{w}^{2} I_{n})$.

As shown in Table~\ref{Table-ex2}, the number of cases with $S_{I} = 0$ among 1000 random systems occurs more frequently than those with $S_{II} = 0$. Moreover, as the noise level increases, the number of cases with $S = 0$ also rises.

\begin{figure}[!ht]
    \centering
    \begin{subfigure}{0.22\textwidth}
        \centering
        \includegraphics[width=\linewidth]{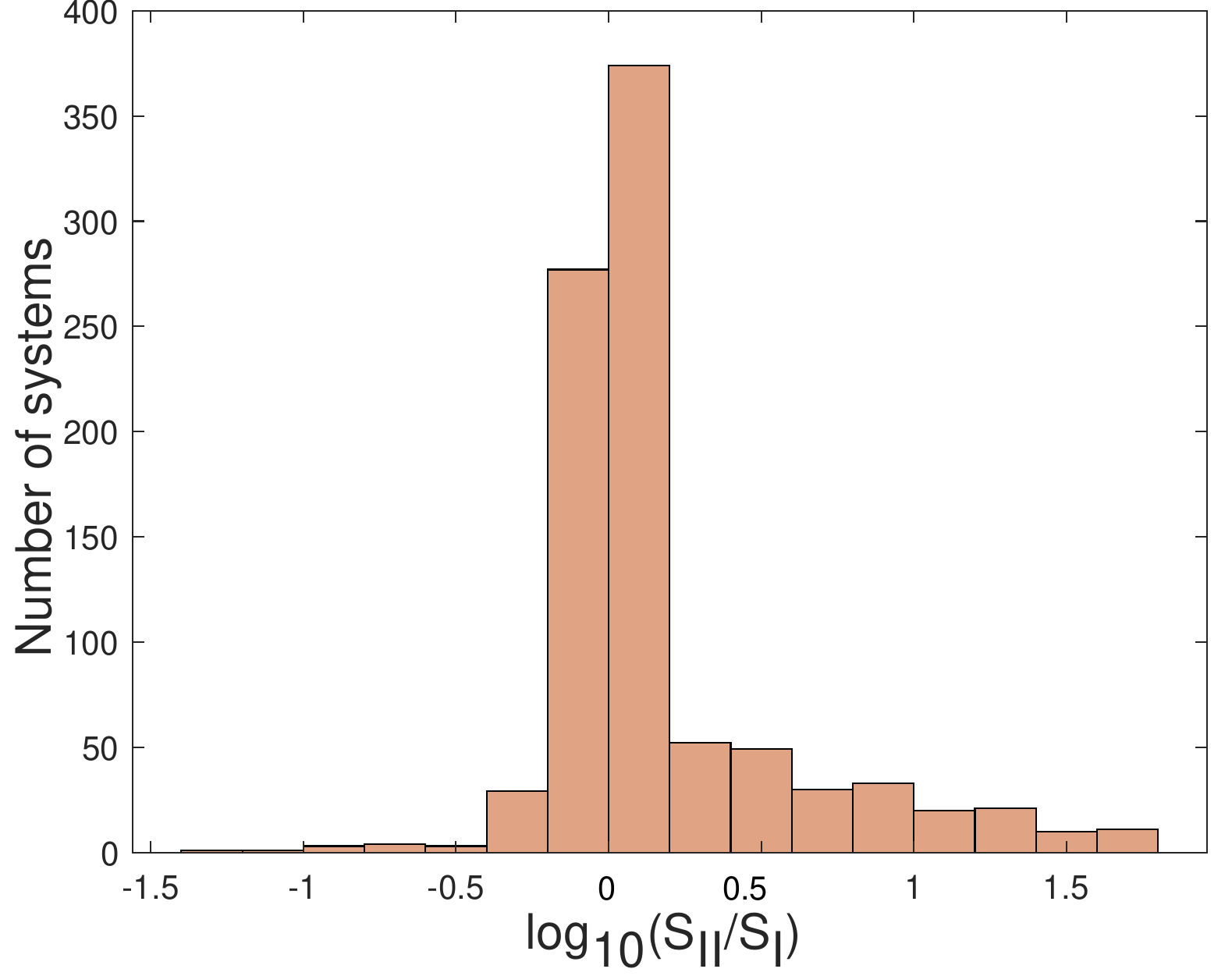}
        \caption{}
    \end{subfigure}
    \quad
    \begin{subfigure}{0.22\textwidth}
        \centering
        \includegraphics[width=\linewidth]{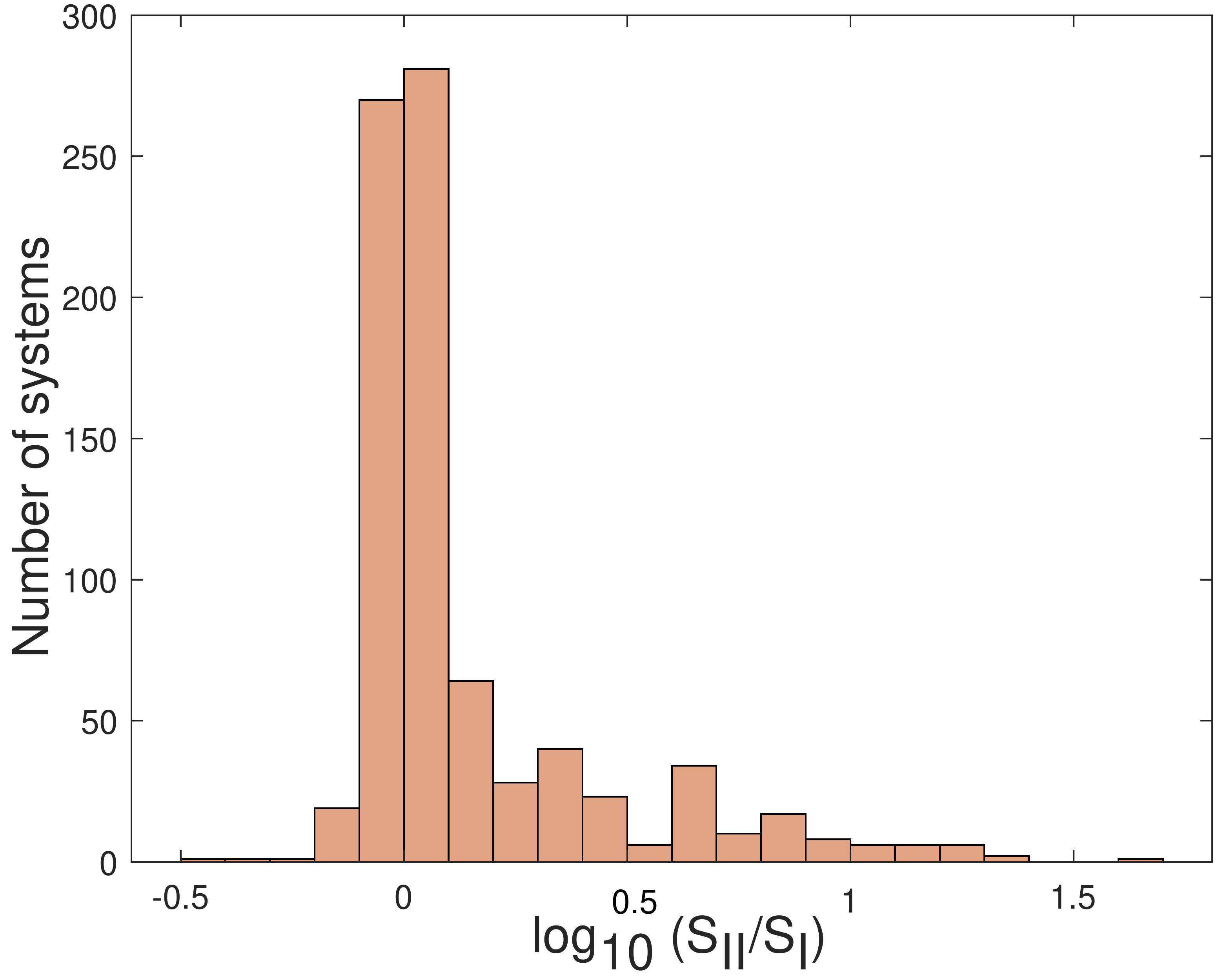}
        \caption{}
    \end{subfigure}
    
\caption{Histograms of $\log_{10}(S_{II}/S_{I})$ for (a) $\sigma_{w} = 0.1$ and (b) $\sigma_{w} = 1$.}
\label{EX4}
\end{figure}

When $S_{I} > 0$ and $S_{II} > 0$, we plot the histograms of $\log_{10}(S_{II}/S_{I})$ for $\sigma_{w} = 0.1$ and $\sigma_{w} = 1$.
As shown in Fig.~\ref{EX4}, the number of cases where $S_{II} > S_{I}$ exceeds those where $S_{II} < S_{I}$, 
demonstrating that our method is more robust in finding a stabilizing controller for random systems.

\begin{table}[htbp]
\centering
\caption{Number of systems (out of 1000) for each case of $S_I$ and $S_{II}$ under different noise levels.}
\label{Table-ex2}
\setlength{\tabcolsep}{12pt}
\renewcommand{\arraystretch}{1.2}
\begin{tabular}{c|c|c}
\Xhline{1pt}
 & $\sigma_{w}=0.1$ & $\sigma_{w}=1.0$ \\
\hline
$S_{I}=0,\,S_{II}=0$     & 16   & 76 \\
$S_{I}=0,\,S_{II}>0$ & 62   & 95 \\
$S_{I}>0,\,S_{II}=0$ & 4    & 5  \\
$S_{I}>0,\,S_{II}>0$ & 918  & 824 \\
\Xhline{1pt}
\end{tabular}
\end{table}

\section{CONCLUSIONS}
\label{5}

In this paper, we proposed a new regularization method for direct data-driven LQR control of unknown LTI systems using the regularized covariance parameterization, which remains effective in handling ill-conditioned cases. We further demonstrated its equivalence to the indirect certainty-equivalence LQR combined with Tikhonov regularization.
Moreover, we extended our method to the design of controllers for unknown nonlinear systems using Koopman linear embedding. Simulation results confirmed the effectiveness and benefits of the proposed method.

In future work, we will develop formal conditions under which our method can effectively address unknown linear and nonlinear systems, informed by the experimental observations in this paper.

\vspace{1em}
\bibliographystyle{ieeetr}
\bibliography{ref_region}

\end{document}